% This is a Plain TeX document

\input amssym.def
\input amssym.tex
\magnification=\magstep1
\footline={\tenrm\hfil --\folio--\hfil}
\baselineskip=18pt  \lineskip=3pt minus 2pt
\lineskiplimit=1pt
%\hsize=15true cm  \vsize=24true cm
\hsize=6.5 true in \vsize=8 true in
\tolerance=1600 \mathsurround=5pt

\def\tomb{\phantom{.}%\hfill
\vrule height.3true cm width.3true cm \par\smallskip\noindent}

\def\trm#1"#2" {\smallskip\parindent=.6true in \itemitem{\bf
#1}{\sl #2}
\hfill\parindent=20pt\smallskip\noindent}
\def\ref
#1.{\mathsurround=0pt${}^{#1}\phantom{|}$\mathsurround=5pt}

\pageno=1  \noindent
\centerline{\bf Nonexistence of Finite-dimensional
Quantizations}
\centerline{\bf of a Noncompact Symplectic Manifold}
\vglue .3in
\centerline{Mark J. Gotay\ref 1). and Hendrik Grundling\ref
2). }
\bigskip
\centerline{\ref 1). Department of Mathematics, University
of Hawai`i,}
\centerline{2565 The Mall, Honolulu, HI 96822 USA}
\centerline{email: gotay@math.hawaii.edu}
\centerline{\ref 2). Department of Pure Mathematics,
University of New South Wales,}
\centerline{ P.O. Box 1, Kensington, NSW 2033, Australia.}
\centerline{ email: hendrik@maths.unsw.edu.au}
\medskip
\centerline{October 21, 1997}
\vglue .2in
\itemitem{{\bf Abstract}}{\sl  We prove that there is no
faithful finite-dimensional representation by skew-hermitian
matrices of a ``basic algebra of observables'' ${\cal B}$ on
a noncompact symplectic manifold $M$. Consequently there
exists no finite-dimensional quantization of {\rm{any}}
Lie subalgebra of the Poisson algebra $C^\infty(M)$
containing ${\cal B}$.}
\vfill\eject

\beginsection 1. Introduction

Let $M$ be a connected noncompact symplectic manifold. On
physical grounds  one expects a quantization of $M$, if it
exists, to be infinite-dimensional. This is what we
rigorously prove here, in the framework of the paper [GGT].
Our precise hypotheses are  spelled out below.

A key ingredient in the quantization process is the choice
of a  {\it basic set of observables} in the Poisson
algebra $C^\infty(M)$.           This is a finite-dimensional
linear subspace ${\cal B}$ of $C^\infty(M)$
such that
\itemitem{(B1)} ({\it Completeness}) the Hamiltonian vector
fields $X_f$, $f\in {\cal B}$, are complete,
\itemitem{(B2)} ({\it Transitivity}) $\{X_f\;\big|\; f\in
{\cal B}\}$ spans $TM$, and
\itemitem{(B3)} ({\it Minimality}) ${\cal B}$ is minimal with
respect to these conditions.

\noindent In addition to these conditions we assume in this
paper that
$\cal B$ forms a Lie algebra under the Poisson bracket. We then
refer to $\cal B$ as a {\it basic algebra}. (Note also that unlike
in [GGT], we do not require here that
$1
\in \cal B.$)

Now fix a
basic algebra
${\cal B}$, and let
${\cal O}$ be any Lie subalgebra of $C^\infty(M)$ containing $1$
and
${\cal B}$. Then by a {\it finite-dimensional
 quantization} of the pair $({\cal O},\,{\cal B})$ we mean a Lie
representation
${\cal Q}$ of $\cal O$ by skew-hermitian matrices on $\Bbb
C^n$ such that
\itemitem{(Q1)} ${\cal Q}(1)=I$,
\itemitem{(Q2)} ${\cal Q}({\cal B})$ is irreducible, and
\itemitem{(Q3)} ${\cal Q}$ is faithful on ${\cal B}$.

We refer the reader to [GGT] for a detailed
discussion of these matters. We remark that in the
infinite-dimensional case there are additional conditions
which must be imposed upon a  quantization.
We also elaborate briefly on (Q3). Although faithfulness is
not usually assumed in the definition of a quantization, it seems
to us a reasonable requirement in that a classical observable can
hardly be regarded as ``basic'' in a physical sense if it is
in the kernel of a quantization map. In this case, it cannot
be obtained in any classical limit from the quantum theory.

\beginsection 2. The Obstruction

Given the definitions above, we state our result:
\trm Theorem." Let $M$ be a noncompact symplectic manifold,
${\cal B}$ a basic algebra on $M$, and
${\cal O}$ any Lie subalgebra of $C^{\infty}(M)$ containing
${\cal B}$. Then there is no finite-dimensional quantization
of $({\cal O},\,{\cal B})$."

As the proof will show, we do not need conditions (Q1) or (Q2) to obtain
the theorem. Moreover, the
subalgebra ${\cal O}$ is irrelevant since the proof depends
only on the Lie theoretical properties of the basic algebra
${\cal B}$ and its action on $M$.

\medskip

\noindent {\bf Proof:} We argue by contradiction. Suppose
there exists a finite-dimensional quantization ${\cal Q}$ of
the basic algebra ${\cal B}$. Since ${\cal Q}({\cal B})$
consists of skew-hermitian matrices, it is  completely
reducible. Since ${\cal Q}$ is faithful, one deduces from [V,
Thm 3.16.3] that ${\cal B}$ is reductive, i.e.
${\cal B}=\frak{s}\oplus\frak{z}$ where $\frak{s}$ is
semisimple and
$\frak{z}$ is the center of ${\cal B}$. We show that
$\frak{z}=\{0\}$. Indeed, by the transitivity condition (B2),
the elements of $\frak{z}$ must be constant but, if these
are nonzero, then $\frak{s}$ alone would serve as a basic algebra,
contradicting the minimality condition (B3). Thus
$\frak{z}=\{0\}$ and ${\cal B}=\frak{s}$ is semisimple.

Let $\it B$ be the connected, simply connected Lie group
with Lie algebra ${\cal B}$. We show that $\it B$ is
noncompact. Let $\frak g$ be the Lie algebra $\{X_f\;\big|\;
f\in {\cal B}\}$. By (B1) the vector fields in $\frak g$
are complete and so by [V, Thm. 2.16.13] this infinitesimal
action of $\frak g$ can be integrated to an action of the
connected, simply connected Lie group $\it G$ with Lie
algebra $\frak g$. Condition (B2) implies that this action
is locally transitive and thus
globally transitive as $M$ is connected. Thus the noncompact
manifold $M$ is a homogeneous space  for $\it G$, and so
$\it G$ must be noncompact as well.  Now ${\cal B}$ is
isomorphic either to $\frak g$ or to a central extension of
$\frak g$ by constants. Since ${\cal B}$ is semisimple, the
latter alternative is impossible. Hence
$\it B$ is isomorphic to
$\it G$ and so is
noncompact.

Now consider a unitary representation $\pi$ of $B$ on $\Bbb C^n$.
Decompose $B$ into a product
$B_1 \times \cdots \times B_K$ of simple groups. Then (at least)
one of these, say
$B_1$, must be noncompact. But it is well-known that a
connected, simple, noncompact Lie group has no nontrivial unitary
representations [BR, Thm. 8.1.2]. Thus $\pi|B_1$ is trivial, i.e.
$\pi(b) = I$ for all $b \in B$. Since every finite-dimensional
quantization
$\cal Q$ of
${\cal B}$ is a derived representation of some unitary
representation $\pi$ of $B$, it follows that ${\cal Q}|{{\cal B}}_1
= 0$, and so $\cal Q$ cannot be faithful.\hfill\tomb

\beginsection 3. Discussion

This theorem is complementary to a recent
result of [GGG] which states that there are no nontrivial
quantizations (finite-dimensional or otherwise) of $(P({\cal B}),\,{\cal
B})$ on a {\it compact\/} symplectic manifold $M$, where
$P({\cal B})$ is the Poisson algebra of polynomials generated
by the basic algebra ${\cal B}$. The proof of that result
leaned heavily on the algebraic structure of
$P({\cal B})$; indeed, when $M$ is compact, it turns out that
${\cal B}$ must be compact semisimple, and such algebras {\it do\/}
have faithful finite-dimensional representations by skew-hermitian
matrices. Thus in the compact case, the obstruction to the
existence of a quantization is Poisson, rather than Lie
theoretical. Combining [GGG] with the
present theorem, we can now assert, roughly speaking, that no
symplectic manifold with a basic algebra has a
finite-dimensional  quantization.

We hope to address the issue of whether there are obstructions
in general to infinite-dimensional quantizations of
noncompact symplectic manifolds in future work. Certainly such
obstructions exist in specific examples, such as $\Bbb R^{2n}$
[GGT] and
$T^*\!S^1$ [GG].
%Moreover, from [GGG] it is known that no compact
%symplectic manifold has a nontrivial infinite-dimensional
%quantization. In effect, then, the case of infinite-dimensional
%quantizations of noncompact symplectic manifolds remains the only
%open problem.
This appears to be a difficult problem, however.

\beginsection Acknowledgments

M.J.G. thanks the University of New South Wales and the U.S.
National  Science Foundation (through grant DMS 96-23083)
for their support while  this research was underway.

H.G would like to thank the Australian Research Council for
their support through a research grant.

\beginsection Bibliography

\item{[BR]} A.O. Barut and R. Ra\c czka [1986] {\sl Theory
of Group Representations and Applications} (Second Edition).
World Scientific, Singapore.
\item{[GGG]} M.J. Gotay, J. Grabowski, and H.B. Grundling
[1997] An Obstruction to Quantizing Compact Symplectic
Manifolds. Preprint dg-ga/9706001.
\item{[GG]} M.J. Gotay and  H.B. Grundling [1997] On Quantizing
$T^*\!S^1$. {\sl Rep. Math. Phys.}, to appear.
\item{[GGT]} M.J. Gotay, H.B. Grundling, G.M. Tuynman [1996]
Obstruction Results in Quantization Theory. {\sl J. Nonlinear
Sci.} {\bf 6}, 469--498.
\item{[V]} V.S. Varadarajan [1984] {\sl Lie Groups, Lie
Algebras and Their Representations.} Springer Verlag, New
York.

\bye